\documentclass[11pt,twoside]{article}
\usepackage{amsmath, amssymb, amsthm}
\usepackage{graphicx} 

\theoremstyle{plain}

\newtheorem{thm}{Theorem}[section]
\newtheorem{prop}[thm]{Proposition}
\newtheorem{lemma}[thm]{Lemma}

\newtheorem{corollary}[thm]{Corollary}

\theoremstyle{definition}

\begin{document}

\def\Z{{\mathbb Z}}\def\N{{\mathbb N}} \def\C{{\mathbb C}}
\def\Q{{\mathbb Q}}\def\R{{\mathbb R}} \def\E{{\mathbb E}}
\def\P{{\mathbb P}}

\let\liff\Longleftrightarrow \let\imply\Rightarrow

\def\Proof{\paragraph{Proof.}}
\def\Remark{\paragraph{Remark.}}
\def\endproof{\hfill$\square$}
\def\noproof{\hfill$\square$}
\def\noi{\noindent}

\let\demph\textbf
\def\PP{{\cal P}}
\def\Sym{{\rm Sym}}
\def\BB{{\cal B}}
\def\deg{{\rm deg}}
\def\MM{{\cal M}}
\def\CC{{\cal C}}
\def\SS{{\cal S}}

\title{\bf{Birman's conjecture for singular\\
braids on closed surfaces}}
 
\author{
\textsc{Luis Paris}}

\date{\today}

\maketitle

\begin{abstract} 
Let $M$ be a closed oriented surface of genus $g\ge 1$, let $B_n(M)$ be the braid group of $M$ on 
$n$ strings, and let $SB_n(M)$ be the corresponding singular braid monoid. 
Our purpose in this paper is to prove that the 
desingularization map $\eta: SB_n(M) \to \Z [B_n(M)]$, introduced in the definition 
of the Vassiliev invariants (for braids on surfaces), is 
injective.
\end{abstract}


%%%%%%%%%%%%%%%%%%%%%%%%%%%%%%%%%%%%%%%%%%%%%%%%%%%%%%%%%%%%%%%%%%%%%%%%%%%%%%%%%%%%%%%%%
\section{Introduction}\label{sect1}

Let $M$ be a surface, and let $\PP= \{ P_1, \dots, P_n\}$ be a collection of $n$ distinct 
punctures in the interior of $M$. Define a {\it braid of $M$ on $n$ strings based at $\PP$} to be 
a $n$-tuple $\beta= (b_1, \dots, b_n)$ of disjoint smooth paths in $M \times [0,1]$, called the 
{\it strings} of $\beta$, such that:

\smallskip
$\bullet$ the projection of $b_i(t)$ on the second coordinate is $t$, for all $t \in [0,1]$ and 
all $i \in \{1, \dots, n\}$;

\smallskip
$\bullet$ $b_i(0)=(P_i,0)$ and $b_i(1)=(P_{\zeta(i)},1)$, where $\zeta$ is a permutation of 
$\{1, \dots, n\}$, for all $i=1, \dots, n$.

\smallskip\noindent
The isotopy classes of braids based at $\PP$ form a group, called the {\it braid group of $M$ on $n$ 
strings based at $\PP$}, denoted by $B_n(M)= B_n(M,\PP)$, and whose multiplication is by 
concatenation. Note that this group does not depend on $\PP$, 
up to isomorphism,
but only on the cardinality 
$n=|\PP|$.

The {\it Artin braid group} $B_n$ is defined to be the braid group on $n$ strings of the plane 
$\E^2$. This group was introduced by Artin in 1926 (see \cite{Art1}, \cite{Art2}), and plays a 
prominent r\^ole in many disciplines. The natural extension to braid groups of topological spaces 
(and, in particular, of surfaces) was introduced by Fox and Neuwirth \cite{FoN} in terms of 
configuration spaces. Presentations for braid groups of closed surfaces have been calculated (see 
\cite{Bir2}, \cite{Sco}, \cite{FaV}, \cite{Gm1}, \cite{Bel2}), these groups are strongly related 
to mapping class groups (see \cite{Bir3}), but very few combinatorial properties of them are 
known. Recently, Irmak, Ivanov, and McCarthy \cite{IIM} have shown that all the automorphisms of 
$B_n(M)$ are geometric (i.e. are induced by diffeomorphisms of $M$), 
provided $M$ is an orientable surface of genus $g\ge 2$ and $n\ge 3$. 
Another important result concerning these groups is a generalization of Markov's theorem which 
relates braids on surfaces to links in 3-dimensional manifolds (see \cite{Sko}). A basic 
reference for surface braid groups is \cite{Bir4}.

Vassiliev invariants, also known as finite type invariants, were first introduced by Vassiliev 
\cite{Va1}, \cite{Va2} for knots, but are now also investigated for other ``knot-type'' classes 
such as links, tangles, Artin braids, or braids on surfaces. 
The general approach is as follows. Given a class $\CC$ of ``knot-like'' objects, one extends $\CC$
to a class $\SS\CC$ of ``singular knot-like'' objects, provided with a filtration $\{\SS_d\CC\}_{d=0}^\infty$,
where $d$ indicates the number of singularities, and with a desingularization map $\eta: \SS\CC \to \Z [\CC]$, where 
$\Z [\CC]$ denotes the free $\Z$-module freely generated by $\CC$.
A {\it Vassiliev invariant 
of order $d$} is then defined to be a homomorphism $v: \Z [\CC] \to A$ of $\Z$-modules which vanishes 
on $\eta( \SS_{d+1} \CC)$.

In the case of braids on surfaces, a {\it singular braid of $M$ on $n$ strings based at $\PP$} is 
a $n$-tuple $\beta= (b_1, \dots, b_n)$ of smooth paths in $M \times [0,1]$, called the {\it 
strings} of $\beta$, such that:

\smallskip
$\bullet$ the projection of $b_i(t)$ on the second coordinate is $t$, for all $i\in \{1, \dots, 
n\}$ and all $t \in [0,1]$;

\smallskip
$\bullet$ $b_i(0)=(P_i,0)$ and $b_i(1)=(P_{\zeta(i)},1)$, where $\zeta$ is a permutation of $\{1, 
\dots, n\}$, for all $i\in \{1, \dots, n\}$;

\smallskip
$\bullet$ the strings of $\beta$ intersect transversely in finitely many double points, called 
{\it singular points} of $\beta$.

\smallskip\noindent
The isotopy classes of singular braids based at $\PP$ form a monoid (and not a group), called 
{\it singular braid monoid of $M$ on $n$ strings based at $\PP$}, and denoted by $SB_n(M)$. It 
obviously contains the braid group $B_n(M)$.

Define the {\it order} of a singular braid to be its number of singular points. Consider a 
singular braid $\beta$ of order $d \ge 1$, and take a singular point $P$ of $\beta$. We can 
slightly modify $\beta$ in a small neighborhood of $P$ in order to suppress the singular point. 
Following this modification, we obtain two singular braids of order $d-1$, denoted by $\beta_+$ 
and $\beta_-$, and called {\it resolutions of $\beta$ at $P$}, as illustrated in Figure~1. Let 
$\Z[B_n(M)]$ denote the braid group algebra of $B_n(M)$. Then we define the {\it desingularization 
map} $\eta: SB_n(M) \to \Z[B_n(M)]$ by induction on the order of a singular braid, setting 
$\eta(\beta)=\beta$ if $\beta$ is a non-singular braid, and $\eta(\beta)= \eta(\beta_+) - 
\eta(\beta_-)$ if $\beta$ is a singular braid of order $d\ge 1$, and $\beta_+, \beta_-$ are the 
resolutions of $\beta$ at some singular point. One can easily verify that $\eta: SB_n(M) \to \Z 
[B_n(M)]$ is a well-defined multiplicative homomorphism.

%%%%%%%%%%%%%%%%%%%%%%%%%%%%%%%%%%%%%%%%%%%%%%%
\begin{figure}[ht]
\begin{center}
\includegraphics[width=10cm]{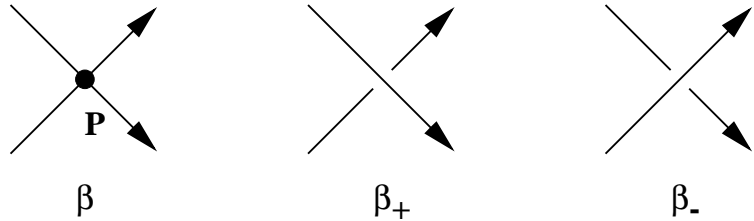}
\end{center}
\caption{Resolutions at a singular point.}
\end{figure}
%%%%%%%%%%%%%%%%%%%%%%%%%%%%%%%%%%%%%%%%%%%%%%%%

Vassilev invariants for braids on closed orientable surfaces of genus $g\ge 1$ have been 
investigated in \cite{GmP}. The main results of \cite{GmP} are:

\smallskip
$\bullet$ the proof that Vassilev invariants separate braids of closed oriented surfaces of 
genus $g\ge 1$;

\smallskip
$\bullet$ the construction of a universal Vassiliev invariant for these braids.

\smallskip\noindent
Note that, by \cite{BeF}, such a universal Vassiliev invariant cannot be functorial (i.e. a 
homomorphism), although, in case of Artin braids, there is a functorial universal invariant (on 
$\Q$) which is defined using the so-called Kontsevich integral (see \cite{Koh}).

A classical question in the subject is to determine, for a given ``knot-like'' class $\CC$, 
whether the desingularization map $\eta: \SS\CC \to \Z [\CC]$ is injective. In the case of Artin 
braids, this question is known as Birman's conjecture (see \cite{Bir1}), and has been recently 
solved in \cite{Par}.

The aim of the present paper is to adapt the techniques of \cite{Par} to braids on closed 
oriented surfaces, and to answer Birman's question in this context. So, our main result is the 
following.

\begin{thm}
Let $M$ be an oriented closed surface of genus $g\ge 1$. Then the desingularization map $\eta: 
SB_n(M) \to \Z [B_n(M)]$ is injective.
\end{thm}

From now on, $M$ denotes a closed oriented surface of genus $g\ge 1$.

One of the keys of the proof of Birman's conjecture in \cite{Par} is that the pure braid group 
(of the plane $\E^2$) can be decomposed as $PB_n= F_{n-1} \rtimes PB_{n-1}$, where $F_{n-1}$ is a 
free group freely generated by some set $\{T_{1\,2}, \dots, T_{1\,n}\}$, and the conjugacy class 
in $F_{n-1}$ of each $T_{1\,j}$ is invariant by the action of $PB_{n-1}$. This fact is not true 
anymore for pure braid groups of closed surfaces. We do have an exact sequence $1 \to R_n(M) \to 
PB_n(M) \to PB_{n-1}(M) \to 1$, where $R_n(M)$ is a free group, but, in general, this exact 
sequence does not split (see \cite{GoG}), and, moreover, the action of $PB_{n-1}(M)$ on the 
abelianization of $R_n(M)$ is not trivial. In order to palliate this difficulty, we replace the 
pure braid group $PB_n(M)$ by the group $K_n(M)$ introduced in \cite{GmP}, and we prove that 
$K_n(M)$ can be decomposed as $K_n(M)= F_n(M) \rtimes K_{n-1}(M)$, where $F_n(M)$ is a free group 
freely generated by some set $\BB'$, and the conjugacy class in $F_n(M)$ of every element of 
$\BB'$ is invariant by the action of $K_{n-1}(M)$. The study of the group $K_n(M)$ is the subject 
of Section 2.

From this point, the proof of Birman's conjecture of \cite{Par} fits quite well to our situation, 
except that we will need to replace the standard homomorphism $\deg: B_n \to \Z$, $\sigma_i \mapsto 1$, 
by some homomorphism $\deg: K_n(M) \to \Z$ (which, by the way, does not extend to $B_n(M)$).

Let $\Gamma$ be a graph, let $X$ be the set of vertices, and let $E=E(\Gamma)$ be the set of edges 
of $\Gamma$. Define the {\it graph monoid} of $\Gamma$ to be the monoid $\MM (\Gamma)$ given by 
the monoid presentation
$$
\MM (\Gamma)= \langle X\ |\ xy=yx \text{ if } \{x,y\} \in E(\Gamma) \rangle^+\,.
$$
In Section 3 we prove that $SB_n(M)$ can be decomposed as $SB_n(M)= \MM (\Omega) \rtimes B_n(M)$, 
where $\MM (\Omega)$ is some graph monoid, and we explain how to use this decomposition to prove 
that $SB_n(M)$ embeds in a group, and to solve the word problem in $SB_n(M)$. In Section 4, we 
show that the desingularization map $\eta: SB_n(M) \to \Z [B_n(M)]$ is injective if and only if 
a certain multiplicative homomorphism $\nu: \MM(\Omega) \to \Z [K_n(M)]$ is injective. Actually, if we 
consider the decomposition $SB_n(M)= \MM (\Omega) \rtimes B_n(M)$ of the previous section, the 
homomorphism $\nu$ turns to be the restriction of $\eta$ to $\MM (\Omega)$. Finally, we prove 
that $\nu: \MM (\Omega) \to \Z [K_n(M)]$ is injective in Section 5.

%%%%%%%%%%%%%%%%%%%%%%%%%%%%%%%%%%%%%%%%%%%%%%%%%%%%%%%%%%%%%%%%%%%%%%%%%%%%%%%%%%%%%%%%%
\section{The braid group and the group $K_n(M)$}\label{sect2}

We represent the surface $M$ by a polygon of $4g$ sides identified as in Figure 2.

%%%%%%%%%%%%%%%%%%%%%%%%%%%%%%%%%%%%%%%%%%%%%%%%%%%%%%%%%%%%%
\begin{figure}[ht]
\begin{center}
\includegraphics[width=6cm]{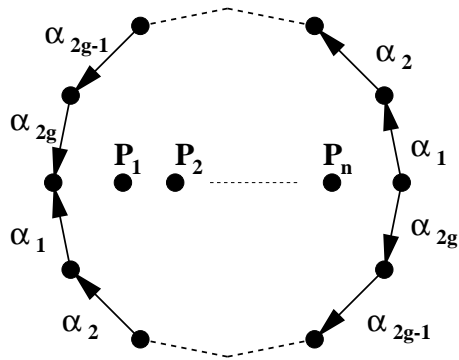}
\end{center}
\caption{The polygon which represents $M$.}
\end{figure}
%%%%%%%%%%%%%%%%%%%%%%%%%%%%%%%%%%%%%%%%%%%%%%%%%%%%%%%%%%%%%

By \cite{Gm1}, $B_n(M)$ has a presentation with generators
$$
\sigma_1, \dots, \sigma_{n-1}, a_1, \dots, a_{2g}\,,
$$
and relations
\[
\begin{array}{lcl}
\text{(R1)}& \sigma_i \sigma_j =\sigma_j \sigma_i&\text{if } |i-j|\ge 2\,,\\
\text{(R2)}& \sigma_i \sigma_{i+1} \sigma_i = \sigma_{i+1} \sigma_i \sigma_{i+1} &\text{if } 1\le 
i\le n-2\,,\\
\text{(R3)}& a_1 \dots a_{2g} a_1^{-1} \dots a_{2g}^{-1} = \sigma_1 \dots \sigma_{n-2} \sigma_{n-
1}^2 \sigma_{n-2} \dots \sigma_1\\
\text{(R4)}& a_r A_{2\,s} = A_{2\,s} a_r &\text{if } 1 \le r,s\le 2g \text{ and } r \neq s\,,\\
\text{(R5)}& (a_1 \dots a_r) A_{2\,r} = \sigma_1^2 A_{2\,r} (a_1 \dots a_r) &\text{if } 1 \le 
r\le 2g\,,\\
\text{(R6)}& a_r \sigma_i = \sigma_i a_r &\text{if } 1\le r\le 2g \text{ and } 2 \le i\le n-
1\,.\\
\end{array}
\]
where
$$
A_{2\,r} = \sigma_1^{-1} (a_1 \dots a_{r-1} a_{r+1}^{-1} \dots a_{2g}^{-1}) \sigma_1^{-1}\,.
$$
We represent the generators of $B_n(M)$ in Figure 3.

%%%%%%%%%%%%%%%%%%%%%%%%%%%%%%%%%%%%%%%%%%%%%%%%%%%%%%%%%%%%%%%%%%%%%%%%%%%%%%
\begin{figure}[ht]
\begin{center}
\includegraphics[width=15cm]{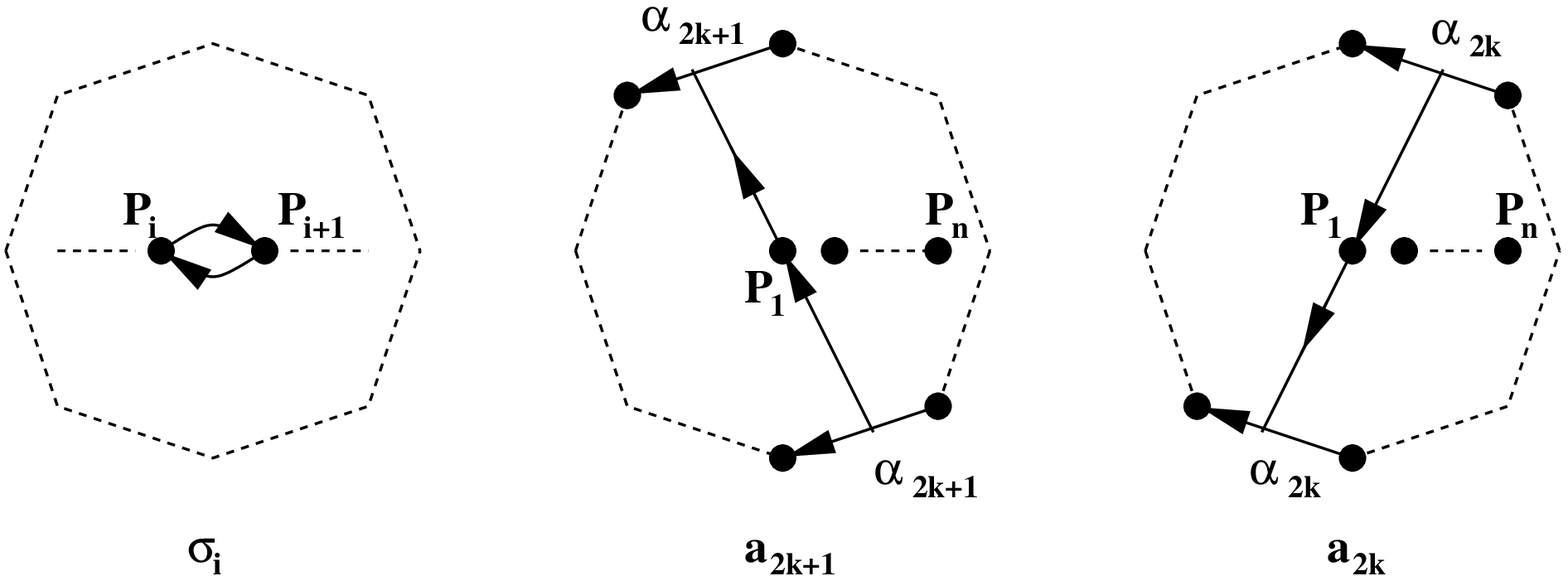}
\end{center}
\caption{Generators of $B_n(M)$.}
\end{figure}
%%%%%%%%%%%%%%%%%%%%%%%%%%%%%%%%%%%%%%%%%%%%%%%%%%%%%%%%%%%%%%%%%%%%%%%%%%%%%%%%

Let $\theta: B_n(M) \to \Sym_n$ be the standard epimorphism defined by $\theta(\sigma_i)= 
(i,i+1)$, for $1 \le i\le n-1$, and $\theta(a_k)=1$ for $1 \le k\le 2g$. The kernel of $\theta$ 
is called the {\it pure braid group of $M$ on $n$ strings based at $\PP$}, and is denoted by 
$PB_n(M)$. Let
\[
\begin{array}{rl}
T_{i\,j}\ =&\sigma_i \dots \sigma_{j-2} \sigma_{j-1}^2 \sigma_{j-2}^{-1} \dots \sigma_i^{-1} 
\quad \text{for }\ 1\le i<j\le n\,,\\
a_{i\,k}\ =& \begin{cases}
\sigma_{i-1}^{-1} \dots \sigma_1^{-1} a_k \sigma_1^{-1} \dots \sigma_{i-1}^{-1} &\text{for } 1 
\le i\le n \text{ and } k \equiv 1\, (\text{mod}\,2)\,,\\
\sigma_{i-1} \dots \sigma_1 a_k \sigma_1 \dots \sigma_{i-1} &\text{for }1\le i\le n \text{ and } 
k \equiv 0\, (\text{mod}\,2)\,.\\
\end{cases}\\
\end{array}
\]
Then $PB_n(M)$ is generated by $\{T_{i\,j}; 1 \le i<j\le n\} \cup \{ a_{i\,k}; 1\le i\le n \text{ 
and } 1 \le k\le 2g\}$ (see \cite{Gm1}). Note that the $T_{i\,j}$'s given in this paper are 
different from the $T_{i\,j}$'s given in \cite{Gm1} and \cite{GmP}, but they are more 
convenient for our purpose. The geometric representations of the generators of $PB_n(M)$ are 
given in Figure 4.

%%%%%%%%%%%%%%%%%%%%%%%%%%%%%%%%%%%%%%%%%%%%%%%%%%%%%%%%%%%%%%%%%%%%%%%%%%ù
\begin{figure}[ht]
\begin{center}
\includegraphics[width=15cm]{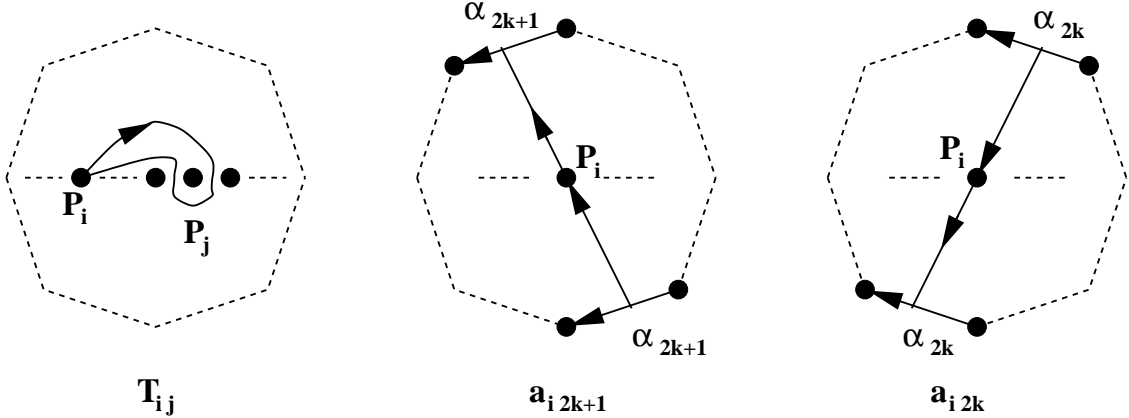}
\end{center}
\caption{Generators of $PB_n(M)$.}
\end{figure}
%%%%%%%%%%%%%%%%%%%%%%%%%%%%%%%%%%%%%%%%%%%%%%%%%%%%%%%%%%%%%%%%%%%%%%%%%%%%

Define the homomorphism $\phi: PB_n(M) \to \pi_1(M)^n$ as follows. Let $\beta= (b_1, \dots, b_n)$ 
be a pure braid. For $i=1, \dots, n$, let $\bar b_i$ be the projection of $b_i$ on the first 
coordinate, and let $\mu_i$ be the element of $\pi_1(M)= \pi_1(M,P_i)$ represented by $\bar b_i$. 
Then $\phi(\beta)= (\mu_1, \dots, \mu_n)$. The kernel of $\phi$ is denoted by $K_n(M)$. By 
\cite{Gol}, this group is the normal closure in $PB_n(M)$ of the subgroup generated by 
$\{T_{i\,j}; 1\le i<j\le n\}$. Let $A_n(M)$ denote the subgroup of $PB_n(M)$ generated by 
$\{a_{i\,k}; 1\le i\le n \text{ and } 1\le k\le 2g\}$. Then $K_n(M)$ is generated by $\{ \alpha 
T_{i\,j} \alpha^{-1} ; \alpha \in A_n(M) \text{ and } 1\le i<j\le n\}$.

Let $\rho: PB_n(M) \to PB_{n-1}(M)$ be the epimorphism which sends a pure braid $\beta= (b_1,b_2, 
\dots, b_n)$ to $(b_2, \dots, b_n)$. Let $\PP_{n-1} = \{P_2, \dots, P_n\}$. By \cite{FN}, the 
kernel of $\rho$ is $\pi_1(M \setminus \PP_{n-1}) = \pi_1(M \setminus \PP_{n-1}, P_1)$. Clearly, 
this is a free group freely generated by $\{T_{1\,j}; 2\le j\le n\} \cup \{a_{1\,k}; 1 \le k\le 
2g\}$.

Let $F_n(M)= K_n(M) \cap \ker \rho = K_n(M) \cap \pi_1(M \setminus \PP_{n-1})$. Then we have the 
following commutative diagram, where all the rows and all the columns are exact.
\[
\begin{array}{ccccccccc}
&&1&&1&&1\\
&&\uparrow&&\uparrow&&\uparrow\\
1&\rightarrow&\pi_1(M)&\rightarrow&\pi_1(M)^n&\rightarrow&\pi_1(M)^{n-1}&\rightarrow&1\\
&&\uparrow&&\uparrow\!\phi&&\uparrow\!\phi\\
1&\rightarrow&\pi_1(M\setminus\PP_{n-1})&\rightarrow&PB_n(M)&\stackrel{\rho}{\rightarrow}&
PB_{n-1}(M)&\rightarrow&1\\
&&\uparrow&&\uparrow&&\uparrow\\
1&\rightarrow&F_n(M)&\rightarrow&K_n(M)&\stackrel{\rho}{\rightarrow}&K_{n-1}(M)&\rightarrow&1\\
&&\uparrow&&\uparrow&&\uparrow\\
&&1&&1&&1\\
\end{array}
\]
Consider the presentation
$$
\pi_1(M) = \langle x_1, \dots, x_{2g}\,|\, x_1x_2 \dots x_{2g}= x_{2g} \dots x_2x_1 \rangle\,.
$$
For all $\gamma \in \pi_1(M)$, we choose a word $\tilde \gamma$ over $\{x_1^{\pm 1}, \dots, 
x_{2g}^{\pm 1}\}$ which represents $\gamma$, and which we call {\it normal form} of $\gamma$. We 
choose the language of normal forms to be prefix-closed, namely, such that any prefix of a normal 
form is also a normal form. For a given word $\omega$ over $\{x_1^{\pm 1}, \dots, x_{2g}^{\pm 
1}\}$, we denote by $\omega_{(i)}$ the word over $\{a_{i\,1}^{\pm 1}, \dots a_{i\,2g}^{\pm 1}\}$ 
obtained from $\omega$ by replacing $x_k^{\pm 1}$ by $a_{i\,k}^{\pm 1}$ for all $k=1, \dots, 2g$.

In spite of the fact that our braids $T_{i\,j}$ are different from the braids $T_{i\,j}$ 
of \cite{GmP}, the proof of the following lemma is exactly the same as the proof of 
\cite{GmP}, Lemma 2.5.

\begin{lemma}\label{Lemma2.1}
For $2\le j\le n$, let
$$
\BB_{1\,j} = \{ \tilde \gamma_{(1)} T_{1\,j} \tilde \gamma_{(1)}^{-1}\ ;\ \gamma \in 
\pi_1(M)\}\,.
$$
Then $F_n(M)$ is a free group freely generated by the disjoint union $\sqcup_{j=2}^n \BB_{1\,j}$. 
\noproof
\end{lemma}

The following lemma is the same as \cite{GmP}, Lemma 2.6. 

\begin{lemma}\label{Lemma2.2}
The epimorphism $\rho: K_n(M) \to K_{n-1}(M)$ admits a section $\iota: K_{n-1}(M) \to K_n(M)$. In 
particular, we can assume that $K_n(M)$ is of the form $K_n(M)= F_n(M) \rtimes K_{n-1}(M)$.
\noproof
\end{lemma}

Now, Proposition 2.3 below is a refinement of \cite{GmP}, Lemma 2.7.

\begin{prop}\label{Proposition2.3}
Let $u \in \sqcup_{j=2}^n \BB_{1\,j}$. Then the conjugacy class of $u$ in $F_n(M)$ is invariant 
by the action of $K_{n-1}(M)$.
\end{prop}

\Proof
Take $u= \tilde \gamma_{(1)} T_{1\,j} \tilde \gamma_{(1)}^{-1} \in \BB_{1\,j}$. Let $r,s \in \{2, 
\dots, n\}$, $r<s$, and let $i\in\{2, \dots, n\}$ and $k \in \{1, \dots, 2g\}$. One can easily 
verify (drawing the corresponding braids) that we have the following relations.
\[
\text{(1)}\quad T_{r\,s} T_{1\,j} T_{r\,s}^{-1} = \begin{cases}
T_{1\,j} &\text{if }1<j<r<s \text{ or } 1<r<s<j\,,\\
T_{1\,j}^{-1} T_{1\,s}^{-1} T_{1\,j} T_{1\,s} T_{1\,j} &\text{if } 1<r=j<s\,,\\
T_{1\,r}^{-1} T_{1\,s}^{-1} T_{1\,r} T_{1\,s} T_{1\,j} T_{1\,s}^{-1} T_{1\,r}^{-1} T_{1\,s} 
T_{1\,r} &\text{if } 1<r<j<s\,,\\
T_{1\,r}^{-1} T_{1\,j} T_{1\,r} &\text{if }1<r<j=s\,.\\
\end{cases}
\]
\[
\text{(2)}\quad a_{i\,k} T_{1\,j} a_{i\,k}^{-1} = \begin{cases}
T_{1,j} &\text{if } i>j\,,\\
T_{1\,j-1} \dots T_{1\,2} (a_{1\,k}^{-1} T_{1\,j} a_{1\,k}) T_{1\,2}^{-1} \dots T_{1\,j-1}^{-1} 
&\text{if } i=j \text{ and } k\equiv 1\,(\text{mod}\,2)\,,\\
a_{1\,k}^{-1} T_{1\,2}^{-1} \dots T_{1\,j-1}^{-1} a_{1\,k} (a_{1\,k}^{-1} T_{1\,j} a_{1\,k})\\ 
\quad
a_{1\,k}^{-1} T_{1\,j-1} \dots T_{1\,2} a_{1\,k} &\text{if } i=j \text{ and } k\equiv 
0\,(\text{mod}\,2)\,,\\
T_{1\,i}^{-1} T_{1\,j} T_{1\,i} &\text{if } 2\le i<j \text{ and } k\equiv1\,(\text{mod}\,2)\,,\\
a_{1\,k}^{-1} T_{1\,2}^{-1} \dots T_{1\,i-1}^{-1} T_{1\,i} T_{1\,i-1} \dots T_{1\,2} a_{1\,k} 
T_{1\,j}\\ 
\quad
a_{1\,k}^{-1} T_{1\,2}^{-1} \dots T_{1\,i-1}^{-1} T_{1\,i}^{-1} T_{1\,i-1} \dots T_{1\,2} 
a_{1\,k} &\text{if } 2\le i<j \text{ and } k\equiv 0\, (\text{mod}\, 2)\,.\\
\end{cases}
\]

Let $z \in \{T_{r\,s}; 2\le r<s\le n\} \cup \{a_{i\,k}; 2\le i \le n \text{ and } 1\le k\le 
2g\}$. Relations (1) and (2) show that there exist $\omega_1 \in F_n(M)$ and $u_1=\tilde 
\mu_{(1)} T_{1\,j} \tilde \mu_{(1)}^{-1} \in \BB_{1\,j}$ such that $z T_{1\,j} z^{-1} = \omega_1 
u_1 \omega_1^{-1}$. let
\[
\begin{array}{c}
\omega_2= z \tilde \gamma_{(1)} z^{-1} \tilde \gamma_{(1)}^{-1} \in F_n(M)\,, \quad \omega_3= 
\tilde \gamma_{(1)} \omega_1 \tilde \gamma_{(1)}^{-1} \in F_n(M)\,,\\
\quad \omega_4=\tilde \gamma_{(1)} \tilde \mu_{(1)} (\widetilde{ \gamma \mu} )_{(1)}^{-1} \in 
F_n(M)\,, \quad u_2= (\widetilde{ \gamma \mu})_{(1)} T_{1\,j} (\widetilde{ \gamma \mu})_{(1)}^{-
1} \in \BB_{1\,j}\,.
\end{array}
\]
Then
$$
zuz^{-1} = (\omega_2 \omega_3 \omega_4) u_2 (\omega_4^{-1} \omega_3^{-1} \omega_2^{-1})\,.
$$

Let $A_{n-1}'(M)$ denote the subgroup of $PB_n(M)$ generated by $\{a_{i\,k}; 2\le i\le n \text{ 
and } 1\le k\le 2g\}$. Let $\alpha \in A_{n-1}'(M)$ and $r,s \in \{2, \dots, n\}$, $r<s$. Then the 
above observations show that there exists $u' \in \BB_{1\,j}$ such that $(\alpha T_{r\,s} 
\alpha^{-1}) u (\alpha T_{r\,s}^{-1} \alpha^{-1})$ is conjugate to $u'$ in $F_n(M)$.

We denote by $H_1(F_n(M))$ the abelianization of $F_n(M)$, and, for $\omega \in F_n(M)$, we 
denote by $[\omega]$ the class of $\omega$ in $H_1(F_n(M))$. Note that $H_1(F_n(M))$ is a free 
abelian group freely generated by $\{ [u]; u\in \sqcup_{j=2}^n \BB_{1\,j}\}$.

Let $r,s \in \{2, \dots, n\}$, $r<s$. Relation (1) shows that there exists $\omega_1 \in F_n(M)$ 
such that $T_{r\,s} T_{1\,j} T_{r\,s}^{-1} = \omega_1 T_{1\,j} \omega_1^{-1}$. Let
$$
\omega_2= T_{r\,s} \tilde \gamma_{(1)} T_{r\,s}^{-1} \tilde \gamma_{(1)}^{-1} \in F_n(M)\,, \quad 
\omega_3= \tilde \gamma_{(1)} \omega_1 \tilde \gamma_{(1)}^{-1} \in F_n(M)\,.
$$
Then
$$
T_{r\,s} u T_{r\,s}^{-1} = (\omega_2 \omega_3) u (\omega_3^{-1} \omega_2^{-1})\,,
$$
hence $[T_{r\,s} u T_{r\,s}^{-1}] = [u]$. This shows that $[T_{r\,s} \omega T_{r\,s}^{-1}] = 
[\omega]$ for all $\omega \in F_n(M)$.

For $\alpha \in A_{n-1}'(M)$ and $\omega \in F_n(M)$, we use the notations $\omega^\alpha= \alpha 
\omega \alpha^{-1}$ and $[\omega]^\alpha = [\alpha \omega \alpha^{-1}]$. Let $\alpha\in A_{n-
1}'(M)$, and $r,s \in \{2, \dots, n\}$, $r<s$. We already know that there exists $u' \in 
\BB_{1\,j}$ such that $(\alpha T_{r\,s} \alpha^{-1}) u (\alpha T_{r\,s}^{-1} \alpha^{-1})$ is 
conjugate to $u'$ in $F_n(M)$. On the other hand,
$$
[u']= [(\alpha T_{r\,s} \alpha^{-1}) u (\alpha T_{r\,s}^{-1} \alpha^{-1})] = [T_{r\,s} 
(u^{\alpha^{-1}}) T_{r\,s}^{-1}]^\alpha = [u^{\alpha^{-1}}]^\alpha = [u]\,,
$$
thus $u'=u$, hence $(\alpha T_{r\,s} \alpha^{-1}) u (\alpha T_{r\,s}^{-1} \alpha^{-1})$ is 
conjugate to $u$ in $F_n(M)$.

Let $\beta \in K_{n-1}(M)$. We choose $\alpha_1, \dots, \alpha_l \in A_{n-1}'(M)$ and $t_1, \dots, 
t_l \in \{T_{r\,s}; 2\le r<s\le n\}$ such that $\beta= \rho((\alpha_1t_1 \alpha_1^{-1}) \dots 
(\alpha_l t_l \alpha_l^{-1}))$. Let
$$
\chi (\beta)= (\alpha_1 t_1 \alpha_1^{-1}) (\alpha_2 t_2 \alpha_2^{-1}) \dots (\alpha_l t_l 
\alpha_l^{-1})\,.
$$
By the above arguments, there exists $\omega_1 \in F_n(M)$ such that $\chi(\beta) u 
\chi(\beta)^{-1} = \omega_1 u \omega_1^{-1}$. Let $\omega_2 =\iota(\beta) \chi(\beta)^{-1} \in 
F_n(M)$. Then $\iota(\beta) u \iota(\beta)^{-1} = (\omega_2 \omega_1) u (\omega_1^{-1} 
\omega_2^{-1})$.
\endproof

\bigskip
Consider the set $\Upsilon=\{ \beta \sigma_i^2 \beta^{-1}; \beta \in B_n(M) \text{ and } 1\le 
i\le n-1\}$, and, for $1\le i<j\le n$, the set $\Upsilon_{i\,j}= \{ \beta T_{i\,j} \beta^{-1}; 
\beta \in PB_n(M)\}$.
Clearly $\Upsilon$ is a generating set of $K_n(M)$, and $\Upsilon_{i\,j}$ is a subset of
$\Upsilon$ for all $1 \le i<j\le n$.

\begin{prop}\label{Proposition2.4}
(1) Let $Z=\oplus_{1\le i<j\le n} \Z e_{i\,j}$
be an (abstract) free abelian group of rank ${n(n-1) \over 2}$, freely generated by some set 
$\{e_{i\,j}; 1\le i<j\le n\}$. There exists a homomorphism $\kappa: K_n(M) \to Z$ which sends $u$ 
to $e_{i\,j}$ for all $u \in \Upsilon_{i\,j}$ and all $1 \le i<j\le n$.

\smallskip
(2) We have the disjoint union $\Upsilon=\sqcup_{i<j} \Upsilon_{i\,j}$.

\smallskip
(3) There exists a homomorphism $\deg : K_n(M) \to \Z$ which sends $u$ to $1$ for all $u \in 
\Upsilon$.
\end{prop}

\Proof
We prove (1) by induction on $n$. Suppose $n=2$. Then $K_2(M)=F_2(M)$ is a free group freely 
generated by $\BB_{1\,2}$, and $\Upsilon= \Upsilon_{1\,2} = \{\alpha u \alpha^{-1}; \alpha\in 
F_2(M) \text{ and }u\in \BB_{1\,2} \}$. Define $\kappa: K_2(M) \to \Z e_{1\,2}$ by $\kappa (u)= 
e_{1\,2}$ for all $u \in \BB_{1\,2}$. Then $\kappa(u)=e_{1\,2}$ for all $u \in \Upsilon_{1\,2}$.

Assume $n>2$. By Proposition 2.3, we have $\Upsilon_{1\,j} = \{ \omega u \omega^{-1}; \omega \in 
F_n(M) \text{ and } u\in \BB_{1\,j}\}$ for all $j=2, \dots, n$. On the other hand, by the 
inductive hypothesis, there exists a homomorphism $\bar \kappa: K_{n-1}(M) \to \oplus_{2\le 
i<j\le n} \Z e_{i\,j}$ which sends $u$ to $e_{i\,j}$ for all $u \in \Upsilon_{i\,j} \subset K_{n-
1}(M)$. Since $K_{n-1}(M)$ acts trivially on the homology of $F_n(M)$, the homomorphism $\bar 
\kappa$ can be extended to a homomorphism $\kappa: K_n(M)= F_n(M) \rtimes K_{n-1}(M) \to Z$ by 
setting $\kappa(u) = e_{1\,j}$ for all $u \in \BB_{1\,j}$ and all $j=2, \dots, n$. Clearly, 
$\kappa$ satisfies $\kappa(u)=e_{i\,j}$ for all $u\in \Upsilon_{i\,j} \subset K_n(M)$ and all 
$1\le i<j\le n$.

Now, we prove (2). The equality $\Upsilon= \cup_{i<j} \Upsilon_{i\,j}$ is obvious. Let $i,j,r,s 
\in \{1, \dots, n\}$, $i<j$, $r<s$, such that $\{i,j\} \neq \{r,s\}$. Let $u \in \Upsilon_{i\,j}$ 
and $v \in \Upsilon_{r\,s}$. Then $\kappa(u)=e_{i\,j} \neq e_{r\,s} =\kappa(v)$, thus $u\neq v$. 
This shows that $\Upsilon_{i\,j} \cap \Upsilon_{r\,s} = \emptyset$.

Consider the homomorphism $\psi: Z \to \Z$ which sends $e_{i\,j}$ to $1$ for all $1\le i<j\le n$. 
Then the homomorphism $\deg= \psi \circ \kappa: K_n(M) \to \Z$ satisfies $\deg(u)=1$ for all 
$u\in \Upsilon$.
\endproof

%%%%%%%%%%%%%%%%%%%%%%%%%%%%%%%%%%%%%%%%%%%%%%%%%%%%%%%%%%%%%%%%%%%%%%%%%%%%%%%%%%%%%%%%%%%
\section{Singular braid monoids}

By \cite{Gm2}, the monoid $SB_n(M)$ has a monoid presentation with generators
$$
\sigma_1^{\pm 1}, \dots, \sigma_{n-1}^{\pm 1}, a_1^{\pm 1}, \dots, a_{2g}^{\pm 1}, \tau_1, \dots, 
\tau_{n-1}\,,
$$
and relations
\[
\begin{array}{lcl}
\text{(R0)}& \sigma_i \sigma_i^{-1} = \sigma_i^{-1} \sigma_i = 1 &\text{if } 1\le i\le n-1\,,\\
&a_k a_k^{-1} = a_k^{-1} a_k = 1 &\text{if } 1\le k\le 2g\,,\\
\text{(R1)-(R6)}& \text{Relations of } B_n(M)\\
\text{(R7)} &\sigma_i \tau_j = \tau_j \sigma_i &\text{if } |i-j|\ge 2\,,\\
\text{(R8)} &\tau_i \tau_j = \tau_j \tau_i &\text{if } |i-j|\ge 2\,,\\
\text{(R9)} &\sigma_i \tau_i = \tau_i \sigma_i &\text{if }1 \le i\le n-1\,,\\
\text{(R10)} &\sigma_i \sigma_j \tau_i = \tau_j \sigma_i \sigma_j &\text{if } |i-j|=1\,,\\
\text{(R11)} &(a_{i\,r} a_{i+1\,r}) \tau_i (a_{i+1\,r}^{-1} a_{i\,r}^{-1}) = \tau_i &\text{if } 
1\le i\le n-1 \text{ and } 1\le r\le 2g\,,\\
\text{(R12)} &\tau_i a_{j\,r} = a_{j\,r} \tau_i &\text{if }j \neq i,i+1 \text{ and } 1\le r\le 
2g\,,\\
\end{array}
\]
where
\[
a_{i\,r} = \begin{cases}
(\sigma_{i-1}^{-1} \dots \sigma_1^{-1}) a_r (\sigma_1^{-1} \dots \sigma_{i-1}^{-1}) &\text{if 
}r\equiv 1\,(\text{mod}\,2)\,,\\
(\sigma_{i-1} \dots \sigma_1) a_r (\sigma_1 \dots \sigma_{i-1}) &\text{if }r\equiv 
0\,(\text{mod}\,2)\,.\\
\end{cases}
\]
The singular braid $\tau_i$ is represented in Figure 5.

%%%%%%%%%%%%%%%%%%%%%%%%%%%%%%%%%%%%%%%%%%%%%%%%%%%%%%%%%%%%%%%%%%%%%%%%%%%%%%%%%%%
\begin{figure}[ht]
\begin{center}
\includegraphics[width=4cm]{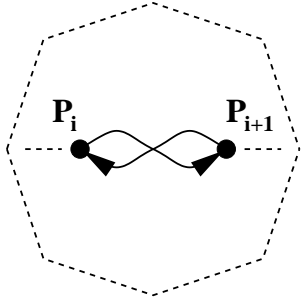}
\end{center}
\caption{The singular braid $\tau_i$.}
\end{figure}
%%%%%%%%%%%%%%%%%%%%%%%%%%%%%%%%%%%%%%%%%%%%%%%%%%%%%%%%%%%%%%%%%%%%%%%%%%%%%%%%%

Now, observe that the desingularization map $\eta: SB_n(M) \to \Z[ B_n(M)]$ is determined by
\[
\begin{array}{cl}
\eta(\sigma_i^{\pm 1}) = \sigma_i^{\pm 1}& \text{for } 1\le i\le n-1\,,\\
\eta( a_k^{\pm 1}) = a_k^{\pm 1} &\text{for } 1\le k\le 2g\,,\\
\eta(\tau_i) = \sigma_i -\sigma_i^{-1} &\text{for } 1\le i\le n-1\,.\\
\end{array}
\]

We shall need for our purpose a slightly different generating set for $SB_n(M)$. Let $\delta_i= 
\sigma_i \tau_i$ for $1\le i\le n-1$. Then $SB_n(M)$ is generated by $\sigma_1^{\pm 1}, \dots, 
\sigma_{n-1}^{\pm 1}, a_1^{\pm 1}, \dots, a_{2g}^{\pm 1}, \delta_1, \dots, \delta_{n-1}$, and has 
a monoid presentation with relations
\[
\begin{array}{lcl}
\text{(R0)}& \sigma_i \sigma_i^{-1} = \sigma_i^{-1} \sigma_i = 1 &\text{if } 1\le i\le n-1\,,\\
&a_k a_k^{-1} = a_k^{-1} a_k = 1 &\text{if } 1\le k\le 2g\,,\\
\text{(R1)-(R6)}& \text{Relations of } B_n(M)\\
\text{(R$7'$)} &\sigma_i \delta_j = \delta_j \sigma_i &\text{if } |i-j|\ge 2\,,\\
\text{(R$8'$)} &\delta_i \delta_j = \delta_j \delta_i &\text{if } |i-j|\ge 2\,,\\
\text{(R$9'$)} &\sigma_i \delta_i = \delta_i \sigma_i &\text{if }1 \le i\le n-1\,,\\
\text{(R$10'$)} &\sigma_i \sigma_j \delta_i = \delta_j \sigma_i \sigma_j &\text{if } |i-j|=1\,,\\
\text{(R$11'$)} &(a_{i\,r} a_{i+1\,r}) \delta_i (a_{i+1\,r}^{-1} a_{i\,r}^{-1}) = \delta_i 
&\text{if } 1\le i\le n-1 \text{ and } 1\le r\le 2g\,,\\
\text{(R$12'$)} &\delta_i a_{j\,r} = a_{j\,r} \delta_i &\text{if }j \neq i,i+1 \text{ and } 1\le 
r\le 2g\,.\\
\end{array}
\]
Moreover, the desingularization map is now determined by
\[
\begin{array}{cl}
\eta(\sigma_i^{\pm 1}) = \sigma_i^{\pm 1}& \text{for } 1\le i\le n-1\,,\\
\eta( a_k^{\pm 1}) = a_k^{\pm 1} &\text{for } 1\le k\le 2g\,,\\
\eta(\delta_i) = \sigma_i^2 -1 &\text{for } 1\le i\le n-1\,.\\
\end{array}
\]

The following proposition is one of the main ingredients of the proof of Theorem 1.1. Its proof can be 
found in \cite{Bel}. It is an extension to surface braid groups of a well-known result on Artin
braid groups due to Fenn, Rolfsen and Zhu \cite{FRZ}.

\begin{prop}
Let $i,j\in \{1, \dots, n-1\}$ and $\beta \in SB_n(M)$. Then the following are equivalent.

\smallskip
(1) $\sigma_i \beta = \beta \sigma_j$;

\smallskip
(2) $\sigma_i^r \beta = \beta \sigma_j^r$ for some $r \in \Z \setminus \{0\}$;

\smallskip
(3) $\tau_i \beta = \beta \tau_j$;

\smallskip
(4) $\tau_i^r \beta = \beta \tau_j^r$ for some $r \in \N \setminus \{0\}$.
\noproof
\end{prop}

\bigskip
Let $\Gamma$ be a graph (with no loop and no multiple edge), let $X$ be the set of vertices, and 
let $E=E(\Gamma)$ be the set of edges of $\Gamma$. Recall that the {\it graph monoid} of $\Gamma$ 
is defined to 
be the monoid $\MM(\Gamma)$ given by the monoid presentation
$$
\MM (\Gamma)= \langle X\ |\ xy=yx \text{ if } \{x,y\} \in E(\Gamma) \rangle^+\,.
$$

The proof of Proposition 3.2 below is exactly the same as the proof of \cite{Par}, Lemma 
2.2.

\begin{prop}
Let $\hat \Omega$ be the graph defined as follows.

\smallskip
$\bullet$ $\hat \Upsilon= \{\alpha \delta_i \alpha^{-1}; \alpha \in B_n(M) \text{ and } 1\le i\le 
n-1\}$ is the set of vertices of $\hat\Omega$;

\smallskip
$\bullet$ $\{\hat u, \hat v\}$ is an edge of $\hat \Omega$ if and only if we have $\hat u \hat v = 
\hat v \hat u$ in $SB_n(M)$.

\smallskip
Let $\Omega$ be the graph defined as follows.

\smallskip
$\bullet$ $\Upsilon= \{ \alpha \sigma_i^2 \alpha^{-1}; \alpha \in B_n(M) \text{ and } 1\le i\le 
n-1\}$ is the set of vertices of $\Omega$;

\smallskip
$\bullet$ $\{u,v\}$ is an edge of $\Omega$ if and only if we have $uv=vu$ in $B_n(M)$.

\smallskip
Then there exists an isomorphism $\varphi: \MM(\hat \Omega) \to \MM (\Omega)$ which sends $\alpha 
\delta_i \alpha^{-1} \in \hat \Upsilon$ to $\alpha \sigma_i^2 \alpha^{-1} \in \Upsilon$ for all 
$\alpha \in B_n(M)$ and $1 \le i\le n-1$.
\noproof
\end{prop}

The proof of the next proposition is the same as the proof of \cite{Par}, Lemma 2.3.

\begin{prop}
We have $SB_n(M)= \MM(\hat \Omega) \rtimes B_n(M)$.
\noproof
\end{prop}

\begin{corollary}
The monoid $SB_n(M)$ embeds in a group.
\end{corollary}

\Proof
Let $G(\hat \Omega)$ be the group given by the presentation
$$
G(\hat \Omega)= \langle \hat \Upsilon\ |\ \hat u \hat v = \hat v \hat u \text{ if } \{ \hat u, 
\hat v\} \in E(\hat \Omega) \rangle\,.
$$
Then $\MM(\hat \Omega)$ embeds in $G(\hat \Omega)$ (see \cite{DK}), thus $SB_n(M)=\MM(\hat 
\Omega) \rtimes B_n(M)$ embeds in $G(\hat \Omega) \rtimes B_n(M)$.
\endproof

\Remark
(1) The result of Corollary 3.4 is due to Bellingeri \cite{Bel} but with a different proof.

(2) Propositions 3.2 and 3.3 together with the solution to the word problem for $B_n(M)$ given in 
\cite{Gm1} can be used to solve the word problem in $SB_n(M)$. Another solution to the word 
problem for $SB_n(M)$ can be found in \cite{Bel}.

%%%%%%%%%%%%%%%%%%%%%%%%%%%%%%%%%%%%%%%%%%%%%%%%%%%%%%%%%%%%%%%%%%%%%%%%%%%%%%%%%%%%%%%%%%%%%
\section{Another theorem}

Recall the subset $\Upsilon=\{ \alpha \sigma_i^2 \alpha^{-1}; \alpha \in B_n(M) \text{ and } 1\le 
i\le n-1\}$ of $K_n(M)$, and the graph $\Omega$ of Proposition 3.2. The goal of the present 
section is to show that Theorem 1.1 is equivalent to the following theorem.

\begin{thm}
Let $\nu: \MM(\Omega) \to \Z [K_n(M)]$ be the multiplicative homomorphism defined by $\nu(u)=u-1$ 
for all $u \in \Upsilon$. Then $\nu$ is injective.
\end{thm}

\paragraph{Proof of ``Theorem 1.1 $\Leftrightarrow$ Theorem 4.1''.}
First, assume that $\eta: SB_n(M) \to \Z [B_n(M)]$ is injective. Recall the isomorphism $\varphi: 
\MM(\hat \Omega) \to \MM (\Omega)$ of Proposition 3.2, and the decomposition $SB_n(M)= \MM (\hat 
\Omega) \rtimes B_n(M)$ of Proposition 3.3, and observe that $\nu= \eta \circ \varphi^{-1}$. So, 
$\nu$ is injective.

Now, assume that $\nu: \MM(\Omega) \to \Z [K_n(M)]$ is injective. Let $L_n(M)= \pi_1(M)^n \rtimes 
\Sym_n$. It is shown in \cite{GmP} that the exact sequence $1 \to K_n(M) \to PB_n(M) 
\stackrel{\phi}{\rightarrow} \pi_1(M)^n \to 1$ extends to an exact sequence
$$
1 \to K_n(M) \to B_n(M) \stackrel{\phi}{\rightarrow} L_n(M) \to 1\,.
$$
Take a set-section $f: L_n(M) \to B_n(M)$ of $\phi$, and consider the isomorphism $\Psi: \Z [B_n(M)] \to \Z 
[K_n(M)] \otimes \Z [L_n(M)]$ of $\Z$-modules defined by
$$
\Psi (\beta)= \beta \cdot (f \circ \phi) (\beta)^{-1} \otimes \phi(\beta)\,,
$$
for $\beta \in B_n(M)$. Recall the homomorphism $\deg: K_n(M) \to \Z$ of Proposition 2.4. For 
$k \in \Z$, let $K_n^{(k)}(M)= \{\beta \in K_n(M); \deg (\beta)=k \}$. We have the decomposition
$$
\Z [K_n(M)] = \bigoplus_{k \in \Z} \Z [K_n^{(k)}(M)]\,,
$$
where $\Z [K_n^{(k)}(M)]$ denotes the free $\Z$-module freely generated by $K_n^{(k)}(M)$. Let $P 
\in \Z [K_n(M)] \otimes \Z [L_n(M)]$. We write $P=\sum_{k \in \Z} P_k$, where $P_k \in \Z 
[K_n^{(k)}(M)] \otimes \Z [L_n(M)]$. Then $P_k$ is called the {\it $k$-th component} of $P$.

Let $\gamma, \gamma' \in SB_n(M)$ such that $\eta(\gamma) = \eta(\gamma')$. We write $\gamma= 
\alpha \beta$ and $\gamma'= \alpha' \beta'$, where $\alpha, \alpha' \in \MM( \hat \Omega)$ and 
$\beta, \beta' \in B_n(M)$ (see Proposition 3.3). Let $\beta_1= \beta (f \circ \phi)(\beta)^{-1} 
\in K_n(M)$, and $\beta_2=\phi(\beta) \in L_n(M)$. Observe that
$$
(\Psi \circ \eta) (\gamma)= (\nu \circ \varphi) (\alpha)\cdot \beta_1 \otimes \beta_2\,.
$$
Let $d=\deg(\beta_1)$. Observe also that the $d$-th component of $(\Psi \circ \eta) (\gamma)$ is 
$\pm \beta_1 \otimes \beta_2$, and, for $k<d$, the $k$-th component of $(\Psi \circ \eta) (\gamma)$ 
is $0$. This shows that $(\Psi \circ \eta) (\gamma)$ determines $\beta_1 \otimes \beta_2$, thus 
$\eta(\gamma)$ determines $\Psi^{-1} (\beta_1 \otimes \beta_2) = \beta$. Since $\eta(\gamma)= 
\eta(\gamma')$, it follows that $\beta= \beta'$.

So, multiplying on the right $\gamma$ and $\gamma'$ by $\beta^{-1}$ if necessary, we can assume 
that $\gamma= \alpha \in \MM( \hat \Omega)$ and $\gamma'= \alpha' \in \MM( \hat \Omega)$. Now, 
observe that
$$
(\nu \circ \varphi) (\gamma)= \eta(\gamma)= \eta( \gamma')= (\nu \circ \varphi) (\gamma')\,,
$$
recall that $\varphi$ is an isomorphism, and recall that $\nu$ is assumed to be injective, hence 
$\gamma=\gamma'$.
\endproof  

%%%%%%%%%%%%%%%%%%%%%%%%%%%%%%%%%%%%%%%%%%%%%%%%%%%%%%%%%%%%%%%%%%%%%%%%%%%%%%%%%%%%%%%%%%%%
\section{Proof of Theorem 4.1}

We start with some results (Lemmas 5.1--5.7) that are preliminary results to the proof of Theorem~4.1.

The following lemma is the same as \cite{Par}, Lemma 3.1.

\begin{lemma}
Let $\Gamma$ be a graph, let $X$ be the set of vertices, and let $E=E(\Gamma)$ be the set of 
edges of $\Gamma$. Let $x_1, \dots, x_l, y_1, \dots, y_l \in X$ and $k\in \{1, \dots, l\}$ such 
that:

\smallskip
$\bullet$ $x_1x_2 \dots x_l= y_1y_2 \dots y_l$ (in $\MM (\Gamma)$);

\smallskip
$\bullet$ $y_k=x_1$ and $y_i \neq x_1$ for all $i=1, \dots, k-1$.

\smallskip\noindent
Then $\{y_i,x_1\} \in E(\Gamma)$ for all $i=1, \dots k-1$.
\noproof
\end{lemma}

Lemma 5.2 below is the same as \cite{Par}, Proposition 4.1.

\begin{lemma}
Let $F(X)$ be a free group freely generated by some set $X$, let $Y= \{ gxg^{-1}; g\in F(X)$ 
and $x\in X\}$, let $F^+(Y)$ be the free monoid freely generated by $Y$, and let $\nu: 
F^+(Y) \to \Z [F(X)]$ be the homomorphism defined by $\nu(y)= y-1$ for all $y \in Y$. Then $\nu$ 
is injective.
\noproof
\end{lemma}

Recall that $\Upsilon= \{\beta \sigma_i^2 \beta^{-1}; \beta \in B_n(M) \text{ and } 1\le i\le n-
1\}$ is the set of vertices of $\Omega$, and that $\Upsilon$ is the disjoint union $\Upsilon= 
\sqcup_{i<j} \Upsilon_{i\,j}$, where $\Upsilon_{i\,j} = \{ \beta T_{i\,j} \beta^{-1}; \beta \in 
PB_n(M)\}$ for $1 \le i<j\le n$ (see Proposition 2.4).

\begin{lemma}
Let $i,j,r,s \in \{1, \dots, n\}$ such that $i<j$, $r<s$, $\{i,j\} \neq \{r,s\}$, and $\{i,j\} 
\cap \{r,s\} \neq \emptyset$. Let $\MM [i,j,r,s]$ be the free monoid freely generated by 
$\Upsilon_{i\,j} \cup \Upsilon_{r\,s}$, and let $\bar \nu: \MM [i,j,r,s] \to \Z [K_n(M)]$ be the 
homomorphism defined by $\bar \nu (u)= u-1$ for all $u \in \Upsilon_{i\,j} \cup \Upsilon_{r\,s}$. 
Then $\bar \nu$ is injective.
\end{lemma}

\Proof
Take $\zeta \in \Sym_n$ such that $\zeta(\{i,j\}) = \{1,n\}$ and $\zeta (\{r,s\}) = \{1,n-1\}$. 
Choose $\beta \in B_n(M)$ such that $\theta (\beta)= \zeta$. Then $\beta \Upsilon_{i\,j} \beta^{-
1} = \Upsilon_{1\,n}$ and $\beta \Upsilon_{r\,s} \beta^{-1} = \Upsilon_{1\,n-1}$. So, up to 
conjugation by $\beta$ if necessary, we can assume that $\{i,j\} = \{1,n\}$ and $\{r,s\} = \{1,n-
1\}$.

Recall the group $F_n(M)$ introduced in Section 2, and, for $j=2, \dots, n$, recall the set 
$\BB_{1\,j} = \{\tilde \gamma_{(1)} T_{1\,j} \tilde \gamma_{(1)}^{-1}; \gamma \in \pi_1(M) \}$. 
Then we have:

\smallskip
$\bullet$ $F_n(M)$ is a free group freely generated by $\BB'=\sqcup_{j=2}^n \BB_{1\,j}$ (see Lemma 
2.1);

\smallskip
$\bullet$ $K_n(M) = F_n(M) \rtimes K_{n-1}(M)$ (see Lemma 2.2);

\smallskip
$\bullet$ $\Upsilon_{1\,j} = \{ \omega t \omega^{-1}; \omega \in F_n(M) \text{ and } t \in 
\BB_{1\,j}\}$, for all $j=2, \dots, n$ (see Proposition 2.3).

\smallskip\noindent
Let $\Upsilon'= \sqcup_{j=2}^n \Upsilon_{1\,j}$, and let $F^+(\Upsilon')$ be the free monoid 
freely generated by $\Upsilon'$. By Lemma 5.2, the homomorphism $\nu': F^+(\Upsilon') \to \Z 
[F_n(M)]$, defined by $\nu'(u)=u-1$ for all $u \in \Upsilon'$, is injective. Recall that $\MM 
[1,n,1,n-1]$ is the free monoid freely generated by $\Upsilon_{1\,n} \cup \Upsilon_{1\,n-
1}$. Then $\MM [1,n,1,n-1] \subset F^+(\Upsilon')$, $\Z [F_n(M)] \subset \Z [K_n(M)]$, and $\bar 
\nu: \MM [1,n,1,n-1] \to \Z [K_n(M)]$ is the restriction of $\nu'$ to $\MM [1,n,1,n-1]$, thus 
$\bar \nu$ is injective.
\endproof

\begin{corollary}
Let $i,j \in \{1, \dots, n\}$ such that $i<j$. Let $\MM [i,j]$ be the free monoid freely 
generated by $\Upsilon_{i\,j}$, and let $\bar \nu: \MM [i,j] \to \Z [K_n(M)]$ be the homomorphism 
defined by $\bar \nu(u) = u-1$ for all $u \in \Upsilon_{i\,j}$. Then $\bar \nu$ is injective.
\noproof
\end{corollary}

The proof of the following lemma is essentially the same as the proof of \cite{Par}, Proposition 
5.1.

\begin{lemma}
Let $F(X)$ be a free group freely generated by some set $X$, let $X_0$ be a subset of $X$, and 
let $\rho: F(X) \to F(X)$ be an automorphism which fixes all the elements of $X_0$, and which 
leaves $F(X \setminus X_0)$ invariant. Let $y_1, \dots, y_l \in \{ \omega x_0 \omega^{-1}; 
\omega\in F(X) \text{ and } x_0 \in X_0\}$. If $\rho(y_1y_2 \dots y_l)= y_1y_2 \dots y_l$, then 
$\rho(y_i)=y_i$ for all $i=1, \dots, l$.
\noproof
\end{lemma}

\begin{lemma}
Assume $n\ge 4$. Let $u_1, \dots, u_l \in \Upsilon_{1\,2}$ and $v \in \Upsilon_{n-1\,n}$. If $v$ 
commutes with $u_1u_2 \dots u_l$ (in $K_n(M)$), then $v$ commutes with $u_i$ for all $i=1, 
\dots, l$.
\end{lemma}

\Proof
Let $\beta_0 \in PB_n(M)$ such that $v= \beta_0 T_{n-1\,n} \beta_0^{-1}$. Up to conjugation of 
$u_1, \dots, u_l$ by $\beta_0^{-1}$ if necessary, we can assume that $\beta_0=1$ and $v=T_{n-
1\,n}$.

Consider the group $F_n(M)$ introduced in Section 2. Recall that:

\smallskip
$\bullet$ $F_n(M)$ is a free group freely generated by the set $\BB'= \sqcup_{j=2}^n \BB_{1\,j}$;

\smallskip
$\bullet$ $K_n(M)= F_n(M) \rtimes K_{n-1}(M)$;

\smallskip
$\bullet$ $\Upsilon_{1\,j} = \{ \omega t \omega^{-1}; \omega \in F_n(M) \text{ and } t \in 
\BB_{1\,j} \}$, for all $j=2, \dots, n$.

\smallskip\noindent
Let $\rho: F_n(M) \to F_n(M)$ be the action of $T_{n-1\,n}$ on $F_n(M)$ by conjugation (namely, 
$\rho (\omega)= T_{n-1\,n} \omega T_{n-1\,n}^{-1}$). Since $n \ge 4$, $T_{n-1\,n}$ commutes with 
$a_k=a_{1\,k}$ for all $k=1, \dots, 2g$. Moreover, we have the following relations (see the proof 
of Proposition 2.3).
\[
T_{n-1\,n} T_{1\,j} T_{n-1\,n}^{-1} = \begin{cases}
T_{1\,j} &\text{if }j<n-1\,,\\
T_{1\,n-1}^{-1} T_{1\,n}^{-1} T_{1\,n-1} T_{1\,n} T_{1\,n-1} &\text{if } j=n-1\,,\\
T_{1\,n-1}^{-1} T_{1\,n} T_{1\,n-1} &\text{if } j=n\,.\\
\end{cases}
\]
Recall also that $\BB_{1\,j} = \{ \tilde \gamma_{(1)} T_{1\,j} \tilde \gamma_{(1)}^{-1}; \gamma 
\in \pi_1(M) \}$ for all $j=2, \dots, n$. These observations imply that $\rho$ fixes all the 
elements of $\BB_{1\,2}$, and leaves invariant the subgroup of $F_n(M)$ generated by $\BB' 
\setminus \BB_{1\,2}$. Now, since $\rho( u_1u_2 \dots u_l)= u_1u_2 \dots u_l$, and $u_1, u_2, 
\dots, u_l \in \Upsilon_{1\,2}= \{ \omega t \omega^{-1}; \omega \in F_n(M) \text{ and } t \in 
\BB_{1\,2}\}$, we conclude by Lemma 5.5 that $\rho(u_i)=u_i$ for all $i=1, \dots, l$.
\endproof

\begin{lemma}
Let $i,j,r,s \in \{ 1, \dots, n\}$ such that $i<j$, $r<s$, and $\{i,j\} \cap \{r,s\} = \emptyset$ 
(in particular, we have $n\ge 4$). Let $\bar \Omega [i,j,r,s]$ be the graph defined as follows.

\smallskip
$\bullet$ $\Upsilon_{i\,j} \cup \Upsilon_{r\,s}$ is the set of vertices of $\bar \Omega 
[i,j,r,s]$;

\smallskip
$\bullet$ $\{u,v\}$ is an edge of $\bar \Omega [i,j,r,s]$ if and only if we have $uv=vu$ in 
$K_n(M)$.

\smallskip\noindent
Let $\MM[i,j,r,s] = \MM (\bar \Omega [i,j,r,s])$, and let $\bar \nu: \MM [i,j,r,s] \to \Z 
[K_n(M)]$ be the homomorphism defined by $\bar \nu(u)=u-1$ for all $u \in \Upsilon_{i\,j} \cup 
\Upsilon_{r\,s}$. Then $\bar\nu$ is injective.
\end{lemma}

\Proof
Take $\zeta \in \Sym_n$ such that $\zeta ( \{i,j\}) = \{n-1,n\}$ and $\zeta(\{r,s\})= \{1,2\}$. 
Choose $\beta \in B_n(M)$ such that $\theta (\beta) = \zeta$. Then $\beta \Upsilon_{i\,j} 
\beta^{-1} = \Upsilon_{n-1\,n}$ and $\beta \Upsilon_{r\,s} \beta^{-1}= \Upsilon_{1\,2}$. So, up 
to conjugation by $\beta$ if necessary, we can assume that $\{i,j\}= \{n-1,n\}$ and 
$\{r,s\}=\{1,2\}$.

Recall the decomposition 
$$
\Z [K_n(M)] = \bigoplus_{k \in \Z} \Z [K_n^{(k)}(M)]\leqno(1)
$$
given in the proof of ``Theorem 1.1 $\Leftrightarrow$ Theorem 4.1'', where $K_n^{(k)}(M)= \{ 
\beta \in K_n(M); \deg(\beta) =k\}$, and $\Z [K_n^{(k)}(M)]$ is the free abelian group freely 
generated by $K_n^{(k)}(M)$. Note that $\deg (u)=1$ for all $u \in \Upsilon$ (see Proposition 
2.4).

Let $\alpha \in \MM [n-1,n,1,2]$. We write $\alpha= u_1u_2 \dots u_l$, where $u_i \in 
\Upsilon_{n-1\,n} \cup \Upsilon_{1\,2}$ for all $i=1, \dots, l$. Define the {\it length} of 
$\alpha$ to be $|\alpha|=l$. We denote by $\bar \alpha$ the element of $K_n(M)$ represented by 
$\alpha$ (i.e. $\bar \alpha= u_1u_2 \dots u_l$ in $K_n(M)$). Let $[1,l]= \{1,2, \dots, l\}$. 
Define a {\it subindex} of $[1,l]$ to be a sequence $I=(i_1,i_2, \dots, i_q)$ such that $i_1,i_2, 
\dots, i_q \in [1,l]$ and $i_1<i_2< \dots <i_q$. The notation $I\prec [1,l]$ means that $I$ is a 
subindex of $[1,l]$. The {\it length} of $I$ is $|I|=q$. For $I=(i_1,i_2, \dots, i_q) \prec 
[1,l]$ we write $\alpha(I)=u_{i_1} u_{i_2} \dots u_{i_q} \in \MM [n-1,n,1,2]$.

Observe that the decomposition of $\bar \nu (\alpha)$ with respect to the direct sum (1) is:
$$
\bar \nu (\alpha)= \sum_{q=0}^l (-1)^{l-q} \sum_{I \prec [1,l],\ |I|=q} \bar 
\alpha(I)\,,\leqno(2)
$$
and
$$
\sum_{I \prec [1,l],\ |I|=q} \bar \alpha(I)\, \in\, \Z [K_n^{(q)}(M)]\,,
$$
for all $q=0,1, \dots, l$.

Let $\alpha'= u_1'u_2' \dots u_k' \in \MM [n-1,n,1,2]$ such that $\bar \nu (\alpha)= \bar \nu 
(\alpha')$. The decomposition given in (2) shows that $k=l$ and
$$
\sum_{I \prec [1,l],\ |I|=q} \bar \alpha(I) = \sum_{I \prec [1,l],\ |I|=q} \bar 
\alpha'(I)\,,\leqno(3)
$$
for all $q=0,1, \dots, l$.

We prove that $\alpha= \alpha'$ by induction on $l$. The cases $l=0$ and $l=1$ being obvious, we 
can assume $l\ge 2$.

Assume first that $u_1=u_1'$. We prove the following equality by induction on $q$.
$$
\sum_{I \prec [2,l],\ |I|=q} \bar \alpha(I) = \sum_{I \prec [2,l],\ |I|=q} \bar \alpha'(I)\,. \leqno(4)
$$
The case $q=0$ being obvious, we can assume $q\ge 1$. Then
\[
\begin{array}{lll}
&{\displaystyle \sum_{I \prec [2,l],\ |I|=q}} \bar \alpha(I)\\
=&{\displaystyle \sum_{I \prec [1,l],\ |I|=q}} \bar \alpha(I)-u_1 \cdot 
{\displaystyle \sum_{I \prec [2,l],\ |I|=q-1}} \bar \alpha(I)\\
=&{\displaystyle \sum_{I \prec [1,l],\ |I|=q}} \bar \alpha'(I)-u_1 \cdot 
{\displaystyle\sum_{I \prec [2,l],\ |I|=q-1}} \bar 
\alpha'(I) &(\text{by induction and (3)})\\
=&{\displaystyle \sum_{I \prec [2,l],\ |I|=q}} \bar \alpha(I)\,.\\
\end{array}
\]
Let $\alpha_1= u_2 \dots u_l$ and $\alpha_1'= u_2' \dots u_l'$. By (4), we have
\[
\begin{array}{rll}
\bar \nu (\alpha_1)\ = &{\displaystyle \sum_{q=0}^{l-1} (-1)^{l-1-q} \sum_{I \prec [2,l],\ |I|=q} \bar 
\alpha(I)}\\
=&{\displaystyle \sum_{q=0}^{l-1} (-1)^{l-1-q} \sum_{I \prec [2,l],\ |I|=q} \bar \alpha'(I)} &=\ \bar \nu 
(\alpha_1')\,.\\
\end{array}
\]
By the inductive hypothesis, it follows that $\alpha_1= \alpha_1'$, hence $\alpha= u_1\alpha_1 = 
u_1 \alpha_1' = \alpha'$.

Now, we consider the general case. Relation (3) applied to $q=1$ gives
$$
\sum_{i=1}^l u_i = \sum_{i=1}^l u_i'\,. \leqno(5)
$$
It follows that there exists a permutation $\zeta \in \Sym_l$ such that $u_i= u_{\zeta (i)}'$ for 
all $i=1, \dots, l$. (Note that the permutation $\zeta \in \Sym_l$ is not necessarily unique. 
Actually, if $u_i=u_j$ for some $i,j\in \{1, \dots l\}$, $i\neq j$, then $\zeta$ is not unique.)

Let $a_1,a_2, \dots, a_p \in [1,l]$, $a_1< a_2< \dots <a_p$, be the indices such that $u_{a_\xi} 
\in \Upsilon_{n-1\,n}$ for all $\xi=1, \dots, p$. Let $I_0=(a_1,a_2, \dots, a_p)$. Recall the 
homomorphism $\kappa: K_n(M) \to Z= \oplus_{i<j} \Z e_{i\,j}$ of Proposition 2.4. Observe that 
$\alpha(I_0) \in \MM [n-1,n]$ and
$$
\bar \nu (\alpha (I_0))= \sum_{k=0}^p (-1)^{p-k} \sum_{\scriptstyle I \prec [1,l],\ |I|=k \atop 
\scriptstyle \kappa (\bar \alpha (I)) \in \Z e_{n-1\,n}} \bar \alpha (I)\,.\leqno(6)
$$
Let $a_1', a_2', \dots, a_p' \in [1,l]$, $a_1'< a_2'< \dots < a_p'$, be the indices such that 
$u_{a_\xi'}' \in \Upsilon_{n-1\,n}$ for all $\xi=1, \dots, p$. Note that $\{ \zeta(a_1'), 
\zeta(a_2'), \dots, \zeta(a_p') \} = \{a_1,a_2, \dots, a_p\}$. Let $I_0'= (a_1', a_2', \dots, 
a_p')$. By (3), we have
$$
\sum_{\scriptstyle I \prec [1,l],\ |I|=k \atop \scriptstyle \kappa (\bar \alpha (I)) \in \Z e_{n-
1\,n}} \bar \alpha (I) = \sum_{\scriptstyle I \prec [1,l],\ |I|=k \atop \scriptstyle \kappa (\bar 
\alpha' (I)) \in \Z e_{n-1\,n}} \bar \alpha' (I)\,, \leqno(7)
$$
for all $k \in \N$, thus, by (6), $\bar \nu (\alpha (I_0)) = \bar \nu (\alpha' (I_0'))$. By 
Corollary 5.4, it follows that $\alpha (I_0) = \alpha' (I_0')$. So, $u_{a_i'}'= u_{a_i}$ for all 
$i=1, \dots, p$, and the permutation $\zeta \in \Sym_l$ can be chosen so that $\zeta(a_i')=a_i$ 
for all $i=1, \dots, p$.

Let $b_1,b_2, \dots, b_q \in [1,l]$, $b_1<b_2< \dots < b_q$, be the indices such that $u_{b_\xi} 
\in \Upsilon_{1\,2}$ for all $\xi=1, \dots, q$. Note that $[1,l]= \{ a_1, \dots, a_p, b_1, \dots, 
b_q\}$. Let $J_0= (b_1,b_2, \dots, b_q)$. Let $b_1', b_2', \dots, b_q' \in [1,l]$, $b_1'< b_2'< 
\dots < b_q'$, be the indices such that $u_{b_\xi'}' \in \Upsilon_{1\,2}$ for all $\xi=1, \dots, 
q$, and let $J_0'= (b_1', b_2', \dots, b_q')$. We also have $\alpha(J_0)= \alpha'(J_0') \in \MM 
[1,2]$, $u_{b_i}= u_{b_i'}'$ for all $i=1, \dots, q$, and $\zeta$ can be chosen so that 
$\zeta(b_i')=b_i$ for all $i=1, \dots, q$.

Without loss of generality, we can assume that $u_1 \in \Upsilon_{n-1\,n}$ (namely, $a_1=1$). Let 
$i \in \{1, \dots, p\}$. We set
\[
\begin{array}{rl}
S(i)\ =&\begin{cases}
0&\text{if } a_i<b_1\,,\\
j&\text{if } b_j<a_i<b_{j+1}\,,\\
q&\text{if }b_q<a_i\,.\\
\end{cases}\\
T(i)\ =&\begin{cases}
0&\text{if } a_i'<b_1'\,,\\
j&\text{if }b_j'<a_i'<b_{j+1}'\,,\\
q&\text{if }b_q'<a_i'\,.\\
\end{cases}\\
\end{array}
\]
Note that $\alpha'= u_{b_1'}' \dots u_{b_{T(1)}'}' u_{a_1'}' \dots = u_{b_1} \dots u_{b_{T(1)}} 
u_{a_1} \dots$. Now, we show that $u_{a_1} =u_1$ commutes with $u_{b_i}$ (in $K_n(M)$ or, 
equivalently, in $\MM [n-1,n,1,2]$) for all $i=1, \dots, T(1)$. It follows that $\alpha' = u_1 
u_{b_1} \dots u_{b_{T(1)}} \dots$, and hence, by the case $u_1=u_1'$ considered before, 
$\alpha=\alpha'$.

Let
\[
\begin{array}{rll}
v_i\ =& u_{b_1} \dots u_{b_{S(i)}} u_{a_i} u_{b_{S(i)}}^{-1} \dots u_{b_1}^{-1} &\in\ 
\Upsilon_{n-1\,n}\,,\\
v_i'\ =& u_{b_1} \dots u_{b_{T(i)}} u_{a_i} u_{b_{T(i)}}^{-1} \dots u_{b_1}^{-1} &\in\ 
\Upsilon_{n-1\,n}\,,\\
\end{array}
\]
for all $i=1, \dots, p$, and let
$$
\gamma= v_1v_2 \dots v_p \in \MM [n-1,n]\,, \quad \gamma'= v_1' v_2' \dots v_p' \in \MM [n-
1,n]\,.
$$
Observe that
\[
\begin{array}{rl}
\bar \nu (\gamma) = &\left( {\displaystyle \sum_{k=0}^p (-1)^{p-k}} {\displaystyle 
\sum_{\scriptstyle I \prec [1,l],\ |I|=k+q \atop \scriptstyle \kappa (\bar \alpha (I)) = ke_{n-
1\,n} +qe_{1\,2}}} \bar \alpha (I) \right) \bar \alpha (J_0)^{-1}\,,\\
\bar \nu (\gamma') = &\left( {\displaystyle \sum_{k=0}^p (-1)^{p-k}} {\displaystyle 
\sum_{\scriptstyle I \prec [1,l],\ |I|=k+q \atop \scriptstyle \kappa (\bar \alpha' (I)) = ke_{n-
1\,n} +qe_{1\,2}}} \bar \alpha' (I) \right) \bar \alpha' (J_0')^{-1}\,.\\
\end{array}
\]
We know that $\alpha (J_0) = \alpha' (J_0')$, and, by (3),
$$
\sum_{\scriptstyle I \prec [1,l],\ |I|=k+q \atop \scriptstyle \kappa (\bar \alpha (I)) = ke_{n-
1\,n} +qe_{1\,2}} \bar \alpha (I) \ = \sum_{\scriptstyle I \prec [1,l],\ |I|=k+q \atop 
\scriptstyle \kappa (\bar \alpha' (I)) = ke_{n-1\,n} +qe_{1\,2}} \bar \alpha' (I)\,,
$$
for all $k=0,1, \dots, p$, thus $\bar \nu (\gamma)= \bar \nu (\gamma')$. By Corollary 5.4, it 
follows that $\gamma=\gamma'$, namely, $v_i=v_i'$ for all $i=1, \dots, p$. So,
$$
u_1= v_1= v_1' = u_{b_1} \dots u_{b_{T(1)}} u_{a_1} u_{b_{T(1)}}^{-1} \dots u_{b_1}^{-1}\,,
$$
thus $u_1$ and $u_{b_1} \dots u_{b_{T(1)}}$ commute (in $K_n(M)$). We conclude by Lemma 5.6 that 
$u_1$ and $u_{b_i}$ commute for all $i=1, \dots, T(1)$.
\endproof

\paragraph{Proof of Theorem 4.1.}
We use the same notation as in the previous proof. Let $\alpha \in \MM (\Omega)$. We write 
$\alpha= u_1u_2 \dots u_l$, where $u_i \in \Upsilon$ for all $i=1, \dots, l$. Observe that
$$
\nu (\alpha)= \sum_{q=0}^l (-1)^{l-q} \sum_{I\prec [1,l],\ |I|=q} \bar \alpha (I)\,,\leqno (1)
$$
and
$$
\sum_{I\prec [1,l],\ |I|=q} \bar \alpha (I)\ \in\ \Z [K_n^{(q)}(M)]\,,
$$
for all $q=0,1, \dots, l$.

Let $\alpha' = u_1' u_2' \dots u_k' \in \MM (\Omega)$ such that $\nu(\alpha) = \nu (\alpha')$. 
Then the decomposition given in (1) shows that $k=l$ and
$$
\sum_{I\prec [1,l],\ |I|=q} \bar \alpha (I)\ = \sum_{I\prec [1,l],\ |I|=q} \bar \alpha' (I)\,, 
\leqno (2)
$$
for all $q=0,1, \dots, l$.

We prove that $\alpha= \alpha'$ by induction on $l$. The cases $l=0$ and $l=1$ being obvious, we 
assume $l \ge 2$.

Assume first that $u_1 = u_1'$. Then, by the same argument as in the previous proof, we have 
$\alpha= \alpha'$.

Now, we consider the general case. Relation (2) applied to $q=1$ gives
$$
\sum_{i=1}^l u_i = \sum_{i=1}^l u_i'\,. \leqno (3)
$$
So, there exists $k \in \{1, \dots, l\}$ such that $u_k'=u_1$ and $u_i' \neq u_1$ for all $i=1, 
\dots, k-1$. We prove that, for $1\le i\le k-1$, $u_i'$ and $u_1=u_k'$ commute (in $K_n(M)$ or, 
equivalently, in $\MM (\Omega)$). It follows that $\alpha'= u_1 u_1' \dots u_{k-1}' u_{k+1}' 
\dots u_l'$, and hence, by the case $u_1=u_1'$ considered before, $\alpha=\alpha'$.

Fix some $t\in \{1, \dots, k-1\}$. Let $i,j,r,s \in \{1, \dots, n\}$ such that $i<j$, $r<s$, 
$u_1=u_k' \in \Upsilon_{i\,j}$, and $u_t' \in \Upsilon_{r\,s}$. There are three possible cases 
that we handle simultaneously:

\smallskip
(1) $\{i,j\} = \{r,s\}$;

\smallskip
(2) $\{i,j\} \neq \{r,s\}$ and $\{i,j\} \cap \{r,s\} \neq \emptyset$;

\smallskip
(3) $\{i,j\} \cap \{r,s\} = \emptyset$.

\smallskip\noindent
Let $\bar \Omega [i,j,r,s]$ be the graph defined as follows.

\smallskip
$\bullet$ $\Upsilon_{i\,j} \cup \Upsilon_{r\,s}$ is the set of vertices of $\bar \Omega 
[i,j,r,s]$;

\smallskip
$\bullet$ $\{u,v\}$ is an edge of $\bar \Omega [i,j,r,s]$ if an only if we have $uv=vu$ in 
$K_n(M)$.

\smallskip\noindent
Let $\MM [i,j,r,s] = \MM (\bar \Omega [i,j,r,s])$, and let $\bar \nu: \MM [i,j,r,s] \to \Z 
[K_n(M)]$ be the homomorphism defined by $\bar \nu (u)=u-1$ for all $u \in \Upsilon_{i\,j} \cup 
\Upsilon_{r\,s}$. Note that, by Corollary 5.4 and Lemma 5.3, $\bar \Omega [i,j,r,s]$ has no edges 
and $\MM [i,j,r,s]$ is a free monoid in Cases (1) and (2). Moreover, the homomorphism $\bar \nu: 
\MM [i,j,r,s] \to \Z [K_n(M)]$ is always injective by Lemmas 5.3 and 5.7 and Corollary 5.4.

Let $a_1, a_2, \dots, a_p \in [1,l]$, $a_1< a_2< \dots <a_p$, be the indices such that $u_{a_\xi} 
\in \Upsilon_{i\,j} \cup \Upsilon_{r\,s}$ for all $\xi=1, \dots, p$. Let $I_0=(a_1, a_2, \dots, 
a_p)$, and let $\alpha (I_0)= u_{a_1} u_{a_2} \dots u_{a_p} \in \MM [i,j,r,s]$. Recall the 
homomorphism $\kappa: K_n(M) \to Z= \oplus_{i<j} \Z e_{i\,j}$ of Proposition 2.4. Observe that
$$
\bar \nu (\alpha (I_0))= \sum_{q=0}^p (-1)^{p-q} \sum_{\scriptstyle I \prec [1,l],\ |I|=q \atop 
\scriptstyle \kappa (\bar \alpha (I)) \in \Z e_{i\,j} + \Z e_{r\,s}} \bar \alpha (I)\,.\leqno(4)
$$
Let $a_1', a_2', \dots, a_p' \in [1,l]$, $a_1'< a_2'< \dots <a_p'$, be the indices such that 
$u_{a_\xi'}' \in \Upsilon_{i\,j} \cup \Upsilon_{r\,s}$ for all $\xi=1, \dots, p$. (Clearly, 
Relation (3) implies that we have as many $a_\xi$'s as $a_\xi'$'s.) Note that $t,k \in \{ a_1', 
\dots, a_p'\}$. Let $I_0'= (a_1', a_2', \dots, a_p')$, and let $\alpha'(I_0') = u_{a_1'}' 
u_{a_2'}' \dots u_{a_p'}' \in \MM [i,j,r,s]$. By (2), we have
$$
\sum_{\scriptstyle I \prec [1,l],\ |I|=q \atop \scriptstyle \kappa (\bar \alpha (I)) \in \Z 
e_{i\,j} + \Z e_{r\,s}} \bar \alpha (I)\ = \sum_{\scriptstyle I \prec [1,l],\ |I|=q \atop 
\scriptstyle \kappa (\bar \alpha' (I)) \in \Z e_{i\,j} + \Z e_{r\,s}} \bar \alpha' (I)\,,
$$
for all $q \in \N$, thus, by (4), $\bar \nu (\alpha (I_0)) = \bar \nu (\alpha' (I_0'))$. Since 
$\bar \nu$ is injective, it follows that $\alpha (I_0)= \alpha' (I_0')$, and we conclude by Lemma 
5.1 that $u_t'$ and $u_k'=u_1$ commute.
\endproof

%%%%%%%%%%%%%%%%%%%%%%%%%%%%%%%%%%%%%%%%%%%%%%%%%%%%%%%%%%%%%%%%%%%%%%%%%%%%%%%%%%%%%%%%

%%%%%%%%%%%%%%%%%%%%%%%%%%%%%%%%%%%%%%%%%%%%%%%%%%%%%%%%%%%%%%%%%%%%%%%%%%%%%%%%%%%%%%%%%
\bigskip\bigskip\noindent
\halign{#\hfill\cr
Luis Paris\cr
Institut de Math\'ematiques de Bourgogne\cr
Universit\'e de Bourgogne\cr
UMR 5584 du CNRS, BP 47870\cr
21078 Dijon cedex\cr
FRANCE\cr
\noalign{\smallskip}
\texttt{lparis@u-bourgogne.fr}\cr}


\begin{thebibliography}{99}

\bibitem{Art1}
E. Artin,
{\it Theory der Z\"opfe},
Abh. Math. Sem. Hamburg {\bf 4} (1926), 47--72.

\bibitem{Art2}
E. Artin,
{\it Theory of braids},
Ann. of Math. {\bf 48} (1947), 101--126.

\bibitem{Bel2}
P. Bellingeri,
{\it On presentation of surface braid groups},
preprint, arXiv: math.GT/0110129.

\bibitem{Bel}
P. Bellingeri,
{\it Centralisers in surface braid groups},
preprint, http://www-fourier.ujf-grenoble.fr/$\sim$bellinge/.

\bibitem{BeF}
P. Bellingeri, L. Funar,
{\it Braids on surfaces and finite type invariants},
preprint, http://www-fourier.ujf-grenoble.fr/$\sim$bellinge/.

\bibitem{Bir2}
J.S. Birman,
{\it On braid groups},
Comm. Pure Appl. Math. {\bf 22} (1968), 41--72.

\bibitem{Bir3}
J.S. Birman,
{\it Mapping class groups and their relationship to braid groups},
Comm. Pure Appl. Math. {\bf 22} (1969), 213--238.

\bibitem{Bir4}
J.S. Birman,
{\it braids, links, and mapping class groups},
Annals of Mathematics Studies, No 82, Princeton University Press, Princeton, N.J., 1974.

\bibitem{Bir1}
J.S. Birman,
{\it New points of view in knot theory},
Bull. Amer. Math. Soc. {\bf 28} (1993), 253--287.

\bibitem{DK}
G. Duchamp, D. Krob,
{\it Free partially commutative structures},
J. Algebra {\bf 156} (1993), 318--361.

\bibitem{FN}
E. Fadell, L. Neuwirth,
{\it Configuration spaces},
Math. Scand. {\bf 10} (1962), 111--118.

\bibitem{FaV}
E. Fadell, J. Van Buskirk,
{\it The braid groups of $E^2$ and $S^2$},
Duke Math. J. {\bf 29} (1962), 243--257.

\bibitem{FRZ}
R. Fenn, D. Rolfsen, J. Zhu,
{\it Centralisers in the braid group and singular braid monoid},
Enseign. Math. {\bf 42} (1996), 75--96.

\bibitem{FoN}
R. Fox, L. Neuwirth,
{\it The braid groups},
Math. Scand. {\bf 10} (1962), 119--126.

\bibitem{Gol}
D.L. Goldsmith,
{\it Homotopy of braids--in answer to a question of E. Artin},
Topology Conference (Virginia Polytech. Inst. and State Univ.,
Blacksburg, Va., 1973), pp. 91--96, Lecture Notes in Math., Vol. 375,
Springer, Berlin, 1974.

\bibitem{GoG}
D.L. Gon\c calves, J. Guaschi,
{\it On the structure of surface pure braid groups},
J. Pure Appl. Algebra {\bf 182} (2003), 33--64.

\bibitem{Gm1}
J. Gonz\'alez-Meneses,
{\it New presentations of surface braid groups},
J. Knot Theory Ramifications {\bf 10} (2001), 431--451.

\bibitem{Gm2}
J. Gonz\'alez-Meneses,
{\it Presentations for the monoids of singular braids on closed surfaces},
Comm. Algebra {\bf 30} (2002), 2829--2836.

\bibitem{GmP}
J. Gonz\'alez-Meneses, L. Paris,
{\it Vassilev invariants for braids on surfaces},
Trans. Amer. Math. Soc., to appear, arXiv: math.GT/0006014.

\bibitem{IIM}
E. Irmak, N.V. Ivanov, J.D. McCarthy,
{\it Automorphisms of surface braid groups},
preprint, arXiv: GT/0306069.

\bibitem{Koh}
T. Kohno,
{\it Vassiliev invariants of braids and iterated integrals},
Arrangements--Tokyo 1998, 157--168,
Adv. Stud. Pure Math., 27,
Kinokuniya, Tokyo, 2000.

\bibitem{Par}
L. Paris,
{\it The proof of Birman's conjecture on singular braid monoids},
preprint, arXiv: math.GR/0306422.

\bibitem{Sco}
G.P. Scott,
{\it Braid groups and the group of homeomorphisms of a surface},
Proc. Cambridge Philos. Soc. {\bf 68} (1970), 605--617.

\bibitem{Sko}
R.K. Skora,
{\it Closed braids in 3-manifolds},
Math. Z. {\bf 211} (1992), 173--187.

\bibitem{Va1}
V.A. Vassilev,
{\it Cohomology of knot spaces},
Theory of singularities and its applications,
23--69, Adv. Soviet Math., 1, Amer. Math. Soc., Providence, RI, 1990.

\bibitem{Va2}
V.A. Vassiliev,
{\it Complements of discriminants of smooth maps: topology and applications},
Translations of Mathematical Monographs, 98, American Mathematical Society,
Providence, RI, 1992.

\end{thebibliography}
\end{document}